\documentclass[a4paper,9pt]{amsart}

\usepackage{graphicx}
\usepackage{natbib}
\usepackage{amsmath}
\usepackage{amssymb}
\usepackage{amsaddr}
\usepackage{booktabs}
\usepackage{mathtools}
\usepackage[T1]{fontenc}
\usepackage{euscript}
\usepackage{tikz}
\usepackage{booktabs}
\usepackage{algorithm}
\usepackage{algorithmicx}
\usepackage{algpseudocode}
\usepackage{amsmath}
\usepackage{enumerate}
\usepackage[hidelinks]{hyperref}

\usetikzlibrary{matrix}

\newcommand\pbatch{\text{$p$-batch, $s_j\leq C$}}
\newtheorem{prop}{Property}
\def\cglb{\text{\tt CG-LB}}
\def\cgub{\text{\tt CG-UB}}

\begin{document}

\title[Column generation for parallel batching]{A solution procedure for
minimizing total completion time in a parallel-batching environment}

\author{A.~Alfieri \and F.~Salassa}
\address{Department of Management and Production Engineering\\
Politecnico di Torino, Italy}
\email{arianna.alfieri@polito.it, fabio.salassa@polito.it}

\author{A.~Druetto \and A.~Grosso}
\address{Department of Computer Science\\
Universit\`a di Torino, Italy}
\email{alessandro.druetto@unito.it, andrea.grosso@unito.it}


\begin{abstract} \,
In many manufacturing processes, batch processing is frequently needed
for capacity reasons. This applies both to parallel and serial
batching. However, while the serial batch processing is largely
studied in the literature, as it is related to the setup issues, the
parallel batch processing is less investigated. In parallel batching,
the manufacturing facility (e.g., ovens for burn-in operations) is
able to accommodate, and process, several parts at the same time and
not to exploit such ability leads to a reduction in the capacity of
the manufacturing facility itself, which will be able to process less
parts per time unit.  In this paper, the scheduling problem in
batch processing environments is considered. Specifically a column
generation algorithm is proposed for parallel batching in both single
machine and parallel machine layouts. Numerical results show that 
the proposed algorithm is able to achieve good solutions in reasonable
computation time due to the use of a new developed lower bound 
much stronger than the literature available lower bounds.
\end{abstract}

\maketitle

\section{Introduction}\label{intro}
In manufacturing system management, capacity is a key factor to have supply matching demand, i.e., to have a system able to produce what is needed to satisfy customer demand. 

Several are the factors negatively impacting the system capacity. The most studied ones are those related to system balancing and to part batching when setup times are present, as severe bottlenecks and/or small batches can largely reduce the system capacity, thus leading to the incapacity of the manufacturing system to timely respond to the market demand (\citet{Cachon}). 

Batches induced by setup times are called serial batches and, although they are very important in manufacturing systems, they are not the only type of batches that can be present in the shop floor. Transfer batches and parallel batches can also be found in manufacturing systems, the first being related to the capacity of the material handling resources and the second, as the serial ones, to the capacity of the machines.

Although both serial and parallel batches are related to and affecting machine capacity, their nature is very different. Serial batches are due the presence of setup times, while parallel batches stem from the ability of machines to accommodate and manufacture several jobs at the same time. They are less studied than serial and transfer batches, because they are less frequent; however, they are not less important. 

Specifically, parallel batches can be found in many manufacturing process where heating operations are necessary, such as in mould manufacturing (\citet{Liu}) and semiconductor industry (\citet{Mason}), or when there are sterilization phases (\citet{Ozturk}), just to cite a few examples.

In all these cases, operations takes a quite long time and the machines usually are batch machines that can accommodate several parts and process them simultaneously, exactly to virtually share the long processing time among all the parts processed at the same time. Each part has an individual size and batch machines (e.g., batch oven for heating treatments or autoclaves for sterilization operations) have a limited capacity; therefore, the number of parts that can be in a single batch is limited.

Due to the limited capacity of the batch machine, and then to the limited number of parts that can be accommodated in it, when several jobs have to be processed on the batch machine, they have to be partitioned in several batches. When batches have been created, their processing has to be scheduled on the machine, and this decision is obviously intertwined with batch creation. Moreover, the two decisions (i.e., how to create batches and how to sequence them on the batch machines) strictly depend on the objective the shop floor manager aims at (e.g., minimizing the number of tardy jobs, minimizing the maximum delay, reducing the total flow time, maximizing the machine utilization, etc.).

In this paper, the above described parallel batch problem is considered. Specifically, given a set of jobs all available at the same time, how to partition them in batches and how to sequence batches on machines is addressed with the objective of minimizing the total completion time. 
 

With respect to the current literature, the problem addressed in the paper is the 
same problem as \cite{IranianiFormiche} 
with the main difference that it extended to the parallel machines case. Following the three field notation of \cite{threefield} 
the problems dealt with is the paper can be refereed to as $1|\pbatch|\sum{C_j}$ and $P_m|\pbatch|\sum{C_j}$, and a column generation algorithm has been developed for their solution.

A fundamental work in the field of parallel batch processor scheduling is the one of \cite{Uzsoy94} where a single batch processing machine problem is studied with regards to makespan and total flow time criteria. In particular, in this work, non identical job sizes are taken into account and complexity results are also provided.

Large part of the literature on parallel batching is devoted to the minimization of the makespan criterion (e.g., \cite{Damodaran06}, \cite{Dupont02} and \cite{Rafiee10}) while the total flow time problem have been less studied (\cite{Jolai98} and \cite{IranianiFormiche}).

The work in \cite{Jolai98} and a modified version of the genetic algorithm presented in \cite{Damodaran06} have been used as benchmark procedures to the hybrid max-min ant system presented in \cite{IranianiFormiche}.
Different objective functions dealing with tardiness and lateness have been also addressed (e.g., \cite{Wang11}, \cite{Malapert12} and \cite{Cabo15}).

The other most recent works on single and parallel machines parallel batching problems are \cite{Beldar18}, \cite{Jia18} and \cite{Ozturk17}.
In \cite{Ozturk17}, the authors address a problem with unit size jobs and maximum completion time objective.
In \cite{Beldar18} and \cite{Jia18}, instead, the total completion time criterion is considered, though \cite{Jia18} consider the weighted version and equal processing time for all the jobs, while in \cite{Beldar18} only the single machine case is tackled with constraints on cardinality and size of jobs batches.

The remainder of the paper is structured as follows.
The column generation approach for $1|\pbatch|\sum C_j$ problem is developed
in Section~\ref{single-mach},
while Section~\ref{parallel-mach} presents the extension of the approach to the parallel machine
case, with special attention of the case with identical parallel machines.
Computational results are reported in Section~\ref{results}.  
Section~\ref{concl} concludes the paper discussing the directions for future research.

\section{Single machine models}
\label{single-mach}
In the remainder of the paper, the following notation is used.
$N=\{1,2,\dots,n\}$ denotes the set of jobs to be
scheduled; for each job $j\in N$, its processing time $p_j$ 
and its size $s_j$, both integers, are given.
The machine has a given integer capacity
denoted by $C$.
When a subset of jobs is packed in a batch $B$, 
$p_B=\max\{p_j\colon j\in B\}$ is used to indicates the batch processing time.
Every batch $B$ is required to have $\sum_{j\in B}s_j\leq C$.
The machine processes the jobs in a batch
sequence $S=(B_1,B_2,\dots,B_t)$, where each job $j$ in the $k$-th
batch $B_k$ shares the batch completion time:
$C_j=C_{B_k}=\sum_{l=1}^k p_{B_l}$\ $\forall j\in B_k$.
The $1|\pbatch|\sum C_j$ problem calls for forming the batches and 
sequencing them in order to minimize $f(S)=\sum_{j\in N}C_j$.

The problem is known to be NP-hard~\citep{Uzsoy94}.
A MILP model for this problem, as given
by~\cite{IranianiFormiche} is the
following, where variable $x_{ij}=1$ iff job $i$ is scheduled in the $j$-th
batch. Variables $C_j$ and $c_i$ represent the completion times of the $j$-th
batch and of job $i$, respectively. Variable $P_j$ represents the processing
time of the $j$-th batch.
\begin{align}
\text{minimize}\quad &\sum_{i = 1}^{n} c_i \label{mip-obj}\\
  \llap{subject to\quad}
& \sum_{j = 1}^{n} x_{ij} = 1 && i = 1, \dots, n \label{mip-1}\\
  & \sum_{i = 1}^{n} s_i x_{ij} \leq C && j = 1, \dots, n \label{mip-2}\\
& P_j \geq x_{ij} p_i && i = 1, \dots, n,\ j = 1, \dots, n \label{mip-3}\\
& C_1 \geq P_1 \label{mip-4}\\ 
& C_j \geq C_{j-1} + P_j && j = 2, \dots, n \label{mip-5}\\ 
& c_i \geq C_j - M \left( 1 - x_{ij} \right) && i = 1, \dots, n,\ j = 1, \dots, n \label{mip-6}\\   
& P_j, C_j, c_i \geq 0 && i = 1, \dots, n,\ j = 1, \dots, n \\
  & x_{ij} \in \left\{ 0, 1 \right\} && i = 1, \dots, n,\ j = 1, \dots, n
  \label{mip-vars}
\end{align}

The total completion time is expressed by~\eqref{mip-obj}. Constraint
set~\eqref{mip-1} ensures that each job is assigned exactly to one batch and,
since all the jobs assigned to a batch cannot exceed the batch capacity, 
constraint set~\eqref{mip-2} has to be defined. Constraint set~\eqref{mip-3} represents the fact that the processing time of a batch is the maximum processing
time of all the contained jobs. 
The completion time for the first batch is simply its processing time since it is
the first to be processed by the machine, as stated in constraint~\eqref{mip-4}.
Constraint set~\eqref{mip-5}, instead, ensures that completion time for all the
other batches is evaluated as the sum of its processing time and the
completion time of the precedent batch. Constraint set~\eqref{mip-6} specifies
that the completion time of a job must be the completion time of the
corresponding batch (the constant $M$ must be very large).

Model \eqref{mip-obj}--\eqref{mip-vars} is known to be very weak. A
state-of-the art solver like \texttt{CPLEX} can waste hours over 15-jobs
instances, with optimality gaps at the root branching node of $100\%$.

For this reason, an alternative model is considered, where a batch
sequence is represented as a path on a graph. Let $\EuScript{B} =
\{B\subseteq N\colon\sum_{j\in B}s_j\leq C\}$ be the set of all the
possible batches. Define a (multi)graph $G(V,A)$ with vertex and arc
sets
\begin{align*}
V&=\{1,2,\dots,n+1\},\\
A&=\left\{ (i,k,B) \colon\ i,k \in V; i<k; B\in\EuScript B;
|B|=(k-i)\right\}.
\end{align*}
Each arc in $A$ is a triple $(i,k,B)$ with head $k$ and
tail $i$ and an associated batch $B$ with $(k-i)$ jobs;
it will represent the batch $B$ scheduled in a batch sequence such that
exactly $n-i+1$ jobs are scheduled
from batch $B$ up to the end of the sequence.
For each arc $(i,k,B)$ a cost is defined as $c_{ikB}=(n - i + 1)p_B$.

\begin{figure}[tbp]
  \begin{minipage}{.5\textwidth}
  \begin{tikzpicture}
    \draw node{$V\quad$};
    \matrix (graph) [right, matrix of math nodes, nodes={rectangle, draw}, ampersand replacement=\&]
            {1\&2\&3\&4\&5\&6\&7\&8\&9\&10\&11\\};
    \draw[->] (graph-1-1) edge[bend left=90] node[above]{$B_1=\{5,6,8\}$} (graph-1-4);
    \draw[->] (graph-1-4) edge[bend right=90] node[below]{$B_2=\{3,4,10\}$} (graph-1-7);
    \draw[->] (graph-1-7) edge[bend left=90] node[above right]{$B_3=\{1,2\}$} (graph-1-9);
    \draw[->] (graph-1-9) edge[bend right=90] node[below right]{$B_4=\{7,9\}$} (graph-1-11);
  \end{tikzpicture}
  \end{minipage}\qquad\qquad
  $\begin{aligned}
    N=&\{1,2,\dots,10\}\\
    S=&(B_1,B_2,B_3,B_4)
  \end{aligned}$
  \[\begin{aligned}
  f(S)=&|B_1|p_{B_1}+\\
       &|B_2|(p_{B_1}+p_{B_2})+\\
       &|B_3|(p_{B_1}+p_{B_2} +p_{B_3})+\\
       &|B_4|(p_{B_1}+p_{B_2} +p_{B_3}+p_{B_4})=\\
       &10p_{B_1}+7p_{B_2}+4p_{B_3}+2p_{B_4}.
    \end{aligned}\]
  \caption{Batch sequence as a path on a graph. $N=B_1\cup B_2\cup B_3\cup B_4$.}\label{fig1}
\end{figure}
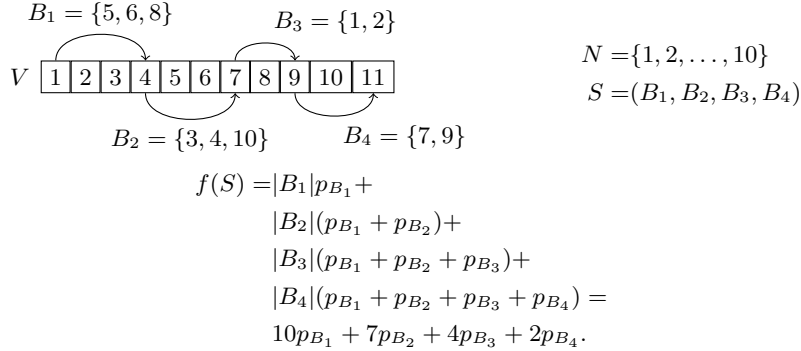

Property \ref{prop1} highlights the relationship between feasible
batches and paths of the above defined graph.
\begin{prop}\label{prop1}
  Each feasible batch sequence $S$ corresponds to a path $P$ in $G(V,A)$ from
  $1$ to $n+1$ such that the set of jobs $N$ is partitioned over
  the arcs in $P$, and $f(S)=\sum\{c_{ikB}\colon (i,k,B)\in P\}$.
\end{prop}
\begin{proof}
  Refer to Figure~\ref{fig1} to fix the idea. A batch sequence
  $S=(B_1,B_2,\dots,B_t)$ is easily seen to be mapped
  (and vice-versa) onto a path
  \[P=[(i_1,k_1,B_1),(i_2,k_2,B_2),\dots,(i_t,k_t,B_t)],\]
  where $B_1,\dots,B_t$ form a partition of $N$,
  $i_1=1$, $k_t=n+1$, $k_\ell=i_\ell+B_\ell$ for $\ell=1,\dots,t$, and 
  $k_\ell=i_{\ell+1}$ for $\ell=1,\dots,t-1$.
  For an arc $(i_q,k_q,B_q)$ in position $q$ on $P$, as claimed above,
  it holds that the number of jobs scheduled from $B_q$ to $B_t$ is 
  \[\sum_{\ell=q}^t|B_\ell|=\sum_{\ell=q}^t(k_\ell-i_\ell)=
  \sum_{\ell=q}^{t-1}i_{\ell+1}+(n+1)-\sum_{\ell=q}^ti_\ell=(n-i_q+1).\]
  The objective function for the batch sequence is
  \begin{align*}
    f(S)&=\sum_{j=1}^nC_j=\sum_{\ell=1}^tC_{B_\ell}|B_\ell|=
    \sum_{\ell=1}^t\Bigg[\Big(\sum_{\kappa=1}^\ell p_{B_\kappa}\Big)|B_\ell|\Bigg]=\\ 
    &=\begin{array}[t]{*{11}l}
    |B_1|p_{B_1}&+\\
    |B_2|p_{B_1}&+&|B_2|p_{B_2}&+&\\
    |B_3|p_{B_1}&+&|B_3|p_{B_2}&+&|B_3|p_{B_3}&+&\\
    \hdotsfor{8}\\
    |B_t|p_{B_1}&+&|B_t|p_{B_2}&+&|B_t|p_{B_3}&+&\dots&+&|B_t|p_{B_t}.
    \end{array}
  \end{align*}
  Adding by column:
  \[\begin{split}
  f(S)=&
  p_{B_1}\sum_{\ell=1}^t|B_\ell|+
  p_{B_2}\sum_{\ell=2}^t|B_\ell|+
  p_{B_3}\sum_{\ell=3}^t|B_\ell|+\dots+p_{B_t}|B_t|=\\
  =&\sum_{q=1}^tp_{B_q}(n-i_q+1)=\sum_{\kern1em\mathclap{(i,k,B)\in P}}c_{ikB}.
  \end{split}
  \]\hfill\qed
\end{proof}

A very large model for the $1|\pbatch|\sum C_j$ problem
includes features of a shortest
path\slash minimum cost flow as well as partition constraints, as follows.
Let $a_B$ be the incidence vector of job set $B$,
$a_{B,j}=\text{($1$ if $j\in B$, else $0$)}$, and
$\mathbf1=(\underbrace{1,1,\dots,1}_n)^{\mathrlap T}$.

\begin{align}
\text{minimize}\quad& \sum_{\mathclap{(i,k,B)\in A}} c_{ikB}x_{ikB} \label{cg:obj}\\
\llap{subject to\quad}
&\sum_{\mathclap{(k,B)\colon\atop(i,k,B)\in A}}x_{ikB}
-\sum_{\mathclap{(k,B)\colon\atop(k,i,B)\in A}}x_{kiB}=
\begin{cases}
    1&i=1\\
    0&i=2,\dots,n\\
    -1&i=n+1
\end{cases}\quad i=1,\dots,n+1\label{cg:flow}\\
&\sum_{\mathclap{(i,k,B)\in A}} a_B x_{ikB} = \mathbf{1}\label{cg:partition} \\
&x_{ikB} \in \left\{ 0, 1 \right\}\qquad (i,k,B)\in A\label{cg:vars}
\end{align}

Variable $x_{ikB}=1$
if arc $(i,k,B)$ is on the path (i.e., batch $B$ is scheduled after $i-1$ jobs);
constraints~\eqref{cg:flow} are flow conservation constraints that
ensure that a unit flow is sent from node $1$ to
node $n+1$, therefore capturing the shortest path problem.
Constraints~\eqref{cg:partition} enforce the requirement that the job set is
exactly partitioned over the arcs selected in the path.

\subsection{Column generation}
The continuous relaxation of~\eqref{cg:obj}--\eqref{cg:vars}, 
where  the integrality constraints~\eqref{cg:vars} are relaxed to
\begin{equation}\label{cg:vars-relaxed}
  x_{ikB}\geq0\qquad(i,k,B)\in A\tag{\ref{cg:vars}$'$},
\end{equation}
is solved by means
of a column generation procedure.
Model~\eqref{cg:obj}--\eqref{cg:vars-relaxed} is the master problem:
a restricted master problem is made of a
subset $A'\subset A$ of arcs.
The dual of~\eqref{cg:obj}--\eqref{cg:vars-relaxed} is
\begin{align}
  \text{maximize}\quad&u_1-u_{n+1}+\sum_{j=1}^n v_j\label{cgd:obj}\\
  \llap{subject to\quad}
  &u_i-u_k+\sum_{j\in B} v_j\leq c_{ikB}
	& (i,k,B)\in A.\label{cgd:dualctr}
\end{align}
Solving the restricted master problem leads to a basic feasible solution
for the master problem and dual variables\slash simplex multipliers
\begin{align*}
  &u_1,\dots,u_{n+1}&\qquad&\text{for constraints~\eqref{cg:flow},}\\
  &v_1,\dots,v_n&\qquad&\text{for constraints~\eqref{cg:partition}.}
\end{align*}
Pricing the arcs $(i,k,B)\in A$ corresponds to finding the most
violated dual constraints~\eqref{cgd:dualctr}.  The strategy developed in this paper is to
price the arcs separately for each pair of indices $i<k$, therefore
determining minimum (possibly negative) reduced costs
\begin{equation}\label{eq:minrcost}
\begin{split}
  \bar c_{ikB^\ast}&=\min_B\Bigg\{
  c_{ikB}-(u_i-u_k)-\sum_{j\in B}v_j\colon \sum_{j\in B}s_j\leq C,\,
  |B|=(k-i)
  \Bigg\}=\\
  &=\min_{B}\Bigg\{p_B(n-i+1)-\sum_{j\in B}v_j\colon \sum_{j\in B}s_j\leq C,\,
  |B|=(k-i)
  \Bigg\}-(u_i-u_k).
\end{split}
\end{equation}

\noindent
Assume to index the jobs in Longest Processing Time (LPT) order,
so that \[p_1\geq p_2\geq\dots\geq p_n.\]
Finding the
batch $B$ that minimizes~\eqref{eq:minrcost} for each given pair
of indices $i<k$ and the batch processing time $p_B$
can be done by exploiting the dynamic programming state space of a family of
cardinality-constrained knapsack problems where items correspond to
jobs. Define, for $r=1,\dots,n$, 
\begin{equation*}
  g_r(\tau,\ell)=\max\Bigg\{\sum_{j=r}^nv_jy_j\colon\sum_{j=r}^ns_jy_j\leq\tau,\
  \sum_{j=r}^ny_j=\ell,\,y_j\in\{0,1\}\Bigg\},
\end{equation*}
where $g_r(\tau,\ell)$ is the optimal value of a knapsack with profits $v_j$ 
and sizes $s_j$, 
limited to items\slash jobs $r,r+1,\dots,n$,
total size $\leq\tau$ and cardinality $\ell$.
Variable $y_j$ is set to 1, i.e., $y_j=1$, iff item\slash job $j$ is included in the solution.

Optimal values for $g_r(\tau, \ell)$ can be recursively
computed as
\[g_r(\tau,\ell)=\max\begin{cases}
g_{r+1}(\tau-s_r,\ell-1)+v_r&(y_r=1)\\
g_{r+1}(\tau,\ell)&(y_r=0)
\end{cases}\]
with boundary conditions
\begin{align*}
  g_r(\tau,1)&=\begin{cases} v_r&\text{if $s_r\leq\tau$\ $(y_r=1)$}\\
  0&\text{otherwise\ $(y_r=0)$}
    \end{cases}&\quad r=1,\dots,n,\,\tau=0,\dots,C\\[2pt]
  g_r(\tau,0)&=0&\quad r=1,\dots,n,\, \tau=0,\dots,C\\
  g_r(\tau,\ell)&=-\infty&\quad \text{if $\ell>n-r+1$ or $\tau<0$.} 
\end{align*}
The corresponding optimal job sets are denoted by  $B_r(\tau,\ell)$;
such sets can be retrieved by backtracking.

\begin{prop}\label{prop2}
  Let $L=\{1\}\cup\{j>1\colon p_j<p_{j-1}\}$.
  For any given pair of indices $i<k$, an arc with minimum reduced cost
  $(i,k,B^\ast)$ is one of
  \begin{equation}\label{eq:redcosts}
    (i,k,B_r(C,k-i))\qquad r\in L.
  \end{equation}
\end{prop}
\begin{proof}
  Every arc $(i,k,B)$ can be shown to have a reduced cost not less than
  some of the arcs in~\eqref{eq:redcosts}.
  Let $\bar c_{ikB}=(n-i+1)p_B-\sum_{j\in B}v_j-(u_i-u_k)$ be the reduced cost of
  an arc $(i,k,B)$.
  Remind that $|B|=k-i$, and the jobs are numbered in non-increasing order of
  processing times.
  Choose $r$ as the smallest job index such that $p_{r}=p_B$. Note that
  $B\subseteq\{r,r+1,\dots,n\}$ and $r\in L$.
  Consider knapsack $g_{r}(C,k-i)$ and the associated
  optimal subset $B_{r}=B_{r}(C,k-i)$.
  The batch $B$ is a feasible solution for 
  knapsack $g_{r}(C,k-i)$, hence $\sum_{j\in B}v_j\leq g_r(C,k-i)$; also, 
  because of the choice of $r$, $p_{B_r}\leq p_r=p_B$. Thus, 
  \[\begin{split}\bar c_{ikB}=&(n-i+1)p_B-\sum_{j\in B}v_j-(u_i-u_k)\geq\\
                         \geq&(n-i+1)p_{B_r}-g_r(C,k-i)-(u_i-u_k)=
                         \bar c_{ikB_r}.
  \end{split}\]
  \hfill\qed
\end{proof}

All the relevant arcs with minimum reduced cost can be generated by
the procedure reported in Algorithm~\ref{algo:pricing} ~(\textsc{NewCols}).
\begin{algorithm}[htbp]
  \caption{Pricing procedure.}\label{algo:pricing}
  \begin{algorithmic}[1]
    \Function{NewCols}{$N$, $C$, $\mathbf u$, $\mathbf v$}
      \Comment{$\mathbf u,\mathbf v=$ vectors of multipliers}
      \State{Sort and renumber
        jobs in $N$ such that $p_1\geq p_2\dots\geq p_n$;}
      \State{Set $H=\emptyset$;}\Comment{Set of negative-reduced cost arcs}
      \For{$\ell=1,\dots,n$}\label{line:for1}
      \State{Set $r:=1$, $\textit{done}:=\textbf{false}$;}
      \While{\textbf{not} \textit{done}}\label{line:while}
        \State{Retrieve $g_r(C,\ell)$ and $B=B_r(C,\ell)$;}
          \label{algo:recursion}
        \For{$i=1,\dots,n-\ell+1$}\label{line:for2}
          \State{Set $k=i+\ell$;}
          \State{Compute $\bar c_{ikB}=p_{B}(n-i+1)-(u_i-u_k)-g_r(C,\ell)$;}
          \label{line:compute-r}
          \If{$\bar c_{ikB}<0$}
            \State{Set $H:=H\cup\{(i,k,B)\}$;}
          \EndIf\label{line:added-arc}
        \EndFor
        \State{Set $r:=\min\{j\colon p_j<p_r\}$;}
        \label{line:update-r}
        \State{If no such index exists, set $\textit{done}:=\textbf{true}$;}
      \EndWhile
    \EndFor
    \State{\textbf{return} $H$;} \EndFunction
  \end{algorithmic}
\end{algorithm}

The size of the state space required for the pricing is bounded by
$\mathcal O(n^2C)$, while the pricing procedure can have two bottlenecks:
\begin{enumerate}[(a)]
\item the $\mathcal O(n^3)$ effort due to the
  three nested loops on lines~\ref{line:for1},
  \ref{line:while}, \ref{line:for2}.
  The \textit{while} on line~\ref{line:while}
  can be executed $n$ times in the worst case;
\item filling the state space $g_r(\tau,\ell)$, which requires at most
  $\mathcal O(n^2C)$ arithmetic operations.  A
  {\em  memoized\/} dynamic programming table is used, so that the execution of the
  top-down recursion for computing an entry $g_r(\tau,\ell)$ is
  deferred until the first time the value is queried.
  Then, the value is kept in storage and accessed in $\mathcal O(1)$ time
  if it is queried again.
\end{enumerate}
Because of these two possible bottlenecks,
the running time of \textsc{NewCols} is bounded from above by
$\mathcal O(\max(n^3, n^2C))$.

\subsection{Heuristic procedure}
The column generation described in the previous section is used to solve the continuous relaxation
of the master problem; once the relaxed optimum has been found, 
the resulting restricted master problem is taken, the variables are set to binary type
and the resulting  MILP is solved by using CPLEX in order to get a heuristic solution
for the master. This is often called ``price and branch'', as
opposed to the exact approach of branch and price.

In order to generate the initial column set,  the jobs are sorted in
shortest processing time (SPT) order and 
all the possible arcs with feasible batches made of SPT-consecutive jobs are generated
(Algorithm~\ref{algo:init}, \textsc{InitCols}).

The complete heuristic procedure is sketched in Algorithm~\ref{algo:heur}.

\begin{algorithm}[tbp]
  \caption{Generation of initial arcs}
  \label{algo:init}
  \begin{algorithmic}
    \Function{InitCols}{$N$, C}
      \State Sort jobs in $N$ such that $p_1\leq p_2\leq\dots\leq p_n$;
      \State Set $H\gets\emptyset$\Comment{Set of initial arcs}
      \For{$j:=1,\dots,n$}
        \State Set $B=\emptyset$;
        \For{$h:=j+1,\dots,n$}
          \If{$\sum_{j\in B}s_j+s_h\leq C$}
            \State Set $B:=B\cup\{h\}$;
            \State Set $H\gets H\cup\{(i,i+|B|,B)\colon i=1,\dots,n-|B|+1\}$;
          \EndIf
        \EndFor
      \EndFor
      \State \textbf{return} $H$;  
    \EndFunction
  \end{algorithmic}
\end{algorithm}

\begin{algorithm}[tbp]
  \caption{Price-and branch procedure}
  \label{algo:heur}
\begin{algorithmic}[1]
\State $A' \gets\textsc{InitCols}(N,C)$;
\State $G(V, A') \gets$ restricted master problem
\While{true}
    \State $z \gets$ continuous optimum of $G(V, A')$
    \State $\mathbf u, \mathbf v\gets$ optimal multipliers
    \State $H \gets\textsc{NewCols}(N,C,\mathbf u, \mathbf v)$
    \If{$|H| = 0$}
        \State $\cglb \gets z$ --- continuous optimum, lower bound
        \State $\cgub \gets$ integer solution computed over $G(V, A')$
        \Comment{Use branch and bound}
        \State \textbf{break}
    \EndIf
    \State $A' \gets A'\cup H$
\EndWhile
\end{algorithmic}
\end{algorithm}

\section{Parallel machine models}
\label{parallel-mach}
Model~\eqref{cg:obj}--\eqref{cg:vars} is readily extended to parallel machine
cases. Consider the fairly general $Rm|\pbatch|\sum C_j$ problem with $m$
parallel unrelated machines.
Let $p_{jh}$ be the processing time of job $j$ on machine $h$.
A special type of arcs with empty batches is added to the graph developed for the single machine case,
using the arc set
\[\begin{split}
  A=&\left\{ (i,k,B) \colon\ i,k \in V; i<k; B\in\EuScript B;
  |B|=(k-i)\right\}\cup\\
  \cup&\left\{(1,k,\emptyset)\colon k=2,\dots,n+1\right\}.
\end{split}\]
Arcs $(i,k,B)\in A$ are given machine-dependent costs
$c_{ikB}^h=p_{Bh}(n-i+1)$, with $p_{Bh}=\max\{p_{jh}\colon j\in B\}$.
Empty arcs $(1,k,\emptyset)$ are given costs
$c_{1k\emptyset}^h=0$, $k=2,\dots,n+1$, $h=1,\dots,m$.
Empty arcs are all added to
the restricted master problem from the beginning, so that they do not need to
be considered in the dynamic programming pricing procedure.  A feasible
solution is made of $m$ batch sequences
\[S_1=(B_1,\dots,B_{t_1}), S_2=(B_{t_1+1},\dots,B_{t_2}),\dots,
  S_m=(B_{t_{m-1}+1},\dots,B_{t_m})\]
processed by the $m$ machines. 
Such batch sequences correspond to $m$ paths (one path for each machine)
$1\to n+1$ on the arcs of which the
set of jobs is exactly partitioned. Such paths will have an empty arc
as first arc. Note that if $(i,k,B)$ is on the $h$-th path,
this means that $n-i+1$ jobs will be scheduled
from $B$ to the end of the $h$-th batch sequence.
Property~\ref{prop1} is easily extended to the multi-machine case. Figure~\ref{fig2} reports  a sketch of the proof with $m=2$. The empty arcs
act as placeholders.

Model~\eqref{cg:obj}--\eqref{cg:vars} can be extended to the parallel machine case
using the multi-commodity features.

\begin{align}
  \text{minimize}\quad& \sum_{h=1}^m\sum_{\,(i,k,B)\mathrlap{\in A}} c_{ikB}^hx_{ikB}^h
  \label{cg2:obj}\\ \noalign{\vskip-1em}
\llap{subject to\quad}
&\sum_{\mathclap{(k,B)\colon\atop(i,k,B)\in A}}x_{ikB}^h
-\sum_{\mathclap{(k,B)\colon\atop(k,i,B)\in A}}x_{kiB}^h=
\begin{cases}
    1&i=1\\
    0&i=2,\dots,n\\
    -1&i=n+1
\end{cases}\quad
\begin{aligned}
  i=&1,\dots,n+1\\
  h=&1,\dots,m
\end{aligned}
\label{cg2:flow}\\
&\sum_{h=1}^m\sum_{(i,k,B)\in A} a_B x_{ikB}^h=\mathbf{1}\label{cg2:partition} \\
&x_{ikB}^h \in \left\{ 0, 1 \right\}\qquad (i,k,B)\in A,\,h=1,\dots,m
\label{cg2:vars}
\end{align}
Here $x_{ikB}^h=1$ iff batch $B$ is on the $h$-th 
path.
Flow conservation constraints~\eqref{cg2:flow}
require that one unit of each commodity is routed
from node $1$ to node $n+1$.
Constraints~\eqref{cg2:partition} enforce the exact partition of the whole
job set across the arcs belonging to the $m$ paths.

In the column generation framework, with each restricted master optimum,
 constraint multipliers are computed:
\begin{align*}
  &u_1^h,u_2^h,\dots,u_{n+1}^h&&\text{for constraints~\eqref{cg2:flow}, $h=1,\dots,m$},\\
  &v_1,v_2,\dots,v_n&&\text{for constraints~\eqref{cg2:partition}.}
\end{align*}
The reduced cost is then separately minimized for each combination of
pair of indices $i<k$ and machine $h$, searching for arcs $(i,k,B)$
with reduced costs 
\[
\bar c_{ikB^\ast}^h=\min_B\Bigg\{p_{Bh}(n-i+1)-\sum_{j\in B}v_j\colon
\sum_{j\in B} s_j\leq C,\,|B|=k-i\Bigg\}-(u_i^h-u_k^h).
\]
This requires calling \textsc{NewCols} $m$ times, once per machine,
since the LPT ordering on each machine is different and so is the state
space $g_r(\tau,\ell)$. Hence, the running time
for pricing raises to $\mathcal O(m\max(n^3,n^2C))$.

A somewhat better situation arises in the case of {\em identical\/}
parallel machines, with problem $Pm|\pbatch|\sum C_j$.
Since each job $j$ has the
same processing time $p_j$ on every machine, the state space
$g_r(\tau,\ell)$ used for pricing is common to all the machines, and a
slightly modified version of \textsc{NewCols} can do  the entire pricing,
still keeping the running time within $\mathcal O(\max(n^3,n^2C))$.
The procedure is reported in Algorithm~\ref{algo:par-pricing}.  The key
observation is that the $p_B$ and $\sum_{j\in B}v_j$ components of the
reduced costs $c_{ikB}^h$ are machine-independent, whereas
only the largest difference $\Delta u_{ik}=\max_h\{(u_i^h-u_k^h)\}$ is
strictly needed in order to compute minimum reduced costs. Such
largest differences are precomputed in time $\mathcal O(mn^2)$ on
line~\ref{getmax}.  For any $(i,k,B)$, let $r$ be the smallest index
such that $p_r=p_B$, and let $B_r=B_r(C,k-i)$;
then, similarly to what proved in
Property~\ref{prop2}:
\[\begin{split}
  \bar c_{ikB}^h=&(n-i+1)p_B-\sum_{j\in B}v_j-(u_i^h-u_k^h)\geq\\
  \geq&(n-i+1)p_{B_r}-g_r(C,k-i)-(u_i^h-u_k^h)\geq\\
  \geq&(n-i+1)p_{B_r}-g_r(C,k-i)-\Delta u_{ik}.
\end{split}
\]


\begin{algorithm}[tbp]
  \caption{Pricing procedure for identical parallel machines.}
  \label{algo:par-pricing}
  \begin{algorithmic}[1]
    \Function{NewCols}{$N$, $C$, $\mathbf u$, $\mathbf v$}
      \Comment{$\mathbf u,\mathbf v=$ vectors of multipliers}
      \State{Sort and renumber jobs in $N$ such
        that $p_1\geq p_2\dots\geq p_n$;}
      \State{Set $\Delta u_{ik}=\max_h\{u_i^h-u_k^h\}$ for $1\leq i<k\leq n+1$;}
      \Comment{$\mathcal O(mn^2)$ time}\label{getmax}
      \State{Set $H=\emptyset$;}\Comment{Set of negative-reduced cost arcs}
      \For{$\ell=1,\dots,n$}\label{line:for1-p}
      \State{Set $r:=1$, $\textit{done}:=\textbf{false}$;}
      \While{\textbf{not} \textit{done}}\label{line:while-p}
        \State{Retrieve $g_r(C,\ell)$ and $B=B_r(C,\ell)$;}
        \For{$i=1,\dots,n-\ell+1$}\label{line:for2-p}
          \State{Set $k=i+\ell$;}
          \State{Compute
            $\bar c_{ikB}=p_{B}(n-i+1)-\Delta u_{ik}-g_r(C,\ell)$;}
          \label{line:compute-r-p}
          \If{$\bar c_{ikB}<0$}
            \State{Set $H:=H\cup\{(i,k,B)\}$;}
          \EndIf\label{line:added-arc-p}
        \EndFor
        \State{Set $r:=\min\{j\colon p_j<p_r\}$;}
        \label{line:update-r-p}
        \State{If no such index exists, set $\textit{done}:=\textbf{true}$;}
      \EndWhile
    \EndFor
    \State{\textbf{return} $H$;}
    \EndFunction
  \end{algorithmic}
\end{algorithm}

Finally note that also taking into account different
capacities for each machines, or even different job sizes on each machine,
simply requires to specialize the knapsack family used. Details are omitted
for the sake of conciseness.

\begin{figure}[tbp]
  \begin{minipage}{.5\textwidth}
  \begin{tikzpicture}
    \draw node{$V\quad$};
    \matrix (graph) [right, matrix of math nodes, nodes={rectangle, draw}, ampersand replacement=\&]
            {1\&2\&3\&4\&5\&6\&7\&8\&9\&10\&11\\};
    \draw[->] (graph-1-1) edge[bend right=90] node[below]{$\emptyset$} (graph-1-5);
    \draw[->] (graph-1-5) edge[bend right=90] node[below]{$B_1$} (graph-1-8);
    \draw[->] (graph-1-8) edge[bend right=90] node[below]{$B_2$} (graph-1-11);
    \draw[->,dashed] (graph-1-1) edge[bend left=90] node[above]{$\emptyset$} (graph-1-7);
    \draw[->, dashed] (graph-1-7) edge[bend left=90] node[above]{$B_3$} (graph-1-9);
    \draw[->, dashed] (graph-1-9) edge[bend left=90] node[above]{$B_4$} (graph-1-11);
  \end{tikzpicture}
  \end{minipage}\qquad\qquad
  $\begin{aligned}
    N=&\{1,2,\dots,10\}\\
    S_1=&(B_1,B_2)\\
    S_2=&(B_3,B_4)\\
    B_1=&\{5,6,8\}\\
    B_2=&\{3,4,10\}\\
    B_3=&\{1,2\}\\
    B_4=&\{7,9\}
  \end{aligned}$
  \[\begin{aligned}
  f(S_1,S_2)=&|B_1|p_{{B_1},1}+\\
             &|B_2|(p_{{B_1},1}+p_{{B_2},1})+\\
             &|B_3|p_{{B_3},2}+\\
             &|B_4|(p_{{B_3},2}+p_{{B_4},2})=\\
             &6p_{B_1,1}+3p_{B_2,1}+4p_{B_3,2}+2p_{B_4,2}.
    \end{aligned}\]
  \caption{Batch sequences on two machines as a 2-path on a graph.
    $N=B_1\cup B_2\cup B_3\cup B_4$.}\label{fig2}
\end{figure}

\section{Computational results}\label{results}

The price and branch heuristics on single machine 
and parallel identical machine instances have been tested on randomly generated instances.
For generating job data, the same approach 
as~\cite{Uzsoy94} and~\cite{IranianiFormiche} have been used. Specifically,
all the job processing times are drawn from
a uniform distribution $p_j \in \left[ 1, 100\right]$,
while job sizes $s_j$ are drawn from four possible uniform distributions,
labeled by $\sigma\in\{\sigma_1,\sigma_2,\sigma_3,\sigma_4\}$:
\[\begin{aligned}
\sigma_1\colon&s_j\in \left[ 1, 10 \right]&\qquad\sigma_3\colon& s_j\in \left[
  3, 10 \right]\\ \sigma_2\colon&s_j\in \left[ 2, 8
  \right]&\qquad\sigma_4\colon&s_j \in \left[ 1, 5 \right].
\end{aligned}\]
In both~\cite{Uzsoy94} and~\cite{IranianiFormiche} the machine capacity
is fixed at $C=10$. Since the pricing procedure has a pseudopolynomial
running time,  instances with $C=30$ and $C=50$ have been also generated in order
to assess how the procedure behaves with a larger capacity.
Single-machine instances have been generated with $n$ ranging from $20$ to $100$ jobs,
and with all four $\sigma$ size distributions. For each $n$, $\sigma$ and~$C$
combinations $10$ random instances have been generated.

With the same job data  the corresponding instances of the
parallel machines problem 
$Pm|\pbatch|\sum C_j$ have been solved for $m=2,3,5$ identical machines. For the parallel machine
case, only the $C=10$ instances have been used.

Both the column-generation based lower bound $\cglb$ and the quality
of the heuristic solution $\cgub$ have been evaluated.
As far as the quality of the lower bound is concerned, the continuous relaxation
of model~\eqref{mip-obj}--\eqref{mip-vars} is not a realistic competitor, zero
being the typical value found by CPLEX at the root branching node.
A more meaningful comparison can be performed against the combinatorial
lower bound proposed by~\cite{Uzsoy94}. Such bound is based on a relaxation of
$1|\pbatch|\sum C_j$ to a preemptive problem on $C$ parallel machines (refer to \cite{Uzsoy94} for details). This lower bound is referred to as \texttt{PR} in the following.

For what the evaluation of $\cgub$ is concerned,
it was difficult to compare the obtained results with the known literature as  neither the test instances nor the
computer codes used by~\cite{Uzsoy94} and~\cite{IranianiFormiche} have been made available.
Hence, some comparison have been made with the results of~\cite{IranianiFormiche},
using instances {\em of the same type,\/} but, for this reason, 
the comparison has to be taken with some care.
On the other hand, when CPLEX is feeded with
model~\eqref{mip-obj}--\eqref{mip-vars} and given some time, its internal
heuristics do generate a number of heuristic solutions, although it has
no chance of certifying optimality. Hence, CPLEX has been run on some set
of instances in order to get heuristic solutions with a time limit of
$300$ seconds.

The times required to compute $\cglb$ and $\cgub$ are separately reported.
The gap between $\cgub$ and $\cglb$ is evaluated as 

\[
gap = \frac{\cgub - \cglb}{\cgub} \cdot 100 \%.
\]

All the tests ran in a Linux environment with \texttt{Intel Core
  i7-6500U CPU @ 2.50GHz} processor; 
\texttt{C++} language has been used for coding the algorithms, and  \texttt{CPLEX} 12.8, 
called directly from \texttt{C++} environment using \texttt{CPLEX}
callable libraries has ben used to solve relaxed and mixed-integer programs.

\subsection{Single machine}
Tables \ref{tab:res10}, \ref{tab:res30} and \ref{tab:res50} show the
results over an increasing number of jobs with batch capacity $C=10$,
$30$ and~$50$, respectively; the $\cgub$was computed using
CPLEX with a time limit of $60$ seconds.  Values are shown as average
over each 10-instance group for the time, and as average, maximum
(worse) and minimum (best) over each 10-instance group for the
gap. Column {\it opt} reports the number of instances (out of 10) in
which the solution can be certified to be the optimum, i.e., in which 
$\cgub=\cglb$.  
The comparison between the $\cglb$ value and 
Uzsoy's \texttt{PR} lower bound is also reported, computing 
the average, maximum and minimum over each 10-instance group of the ratio $\cglb / \texttt{PR}$.

\begin{table}[tbp]
\centering
\begin{tabular}{ll|rr|rrr|rrr|c}
\toprule
\multicolumn{2}{c}{Param} & \multicolumn{2}{c}{Times (s)} & \multicolumn{3}{c}{Gap (\%)} & \multicolumn{3}{c}{$\cglb / \texttt{PR}$} \\
\cmidrule(r){1-2}
\cmidrule(lr){3-4}
\cmidrule(lr){5-7}
\cmidrule(lr){8-10}
$n$ & $\sigma$ & \cglb & \cgub & avg & worst & best & avg & min & max & opt \\
\midrule
20 & $\sigma_1$ & 0.01 & 0.04 & 1.30 & 3.20 & 0.00 & 1.25 & 1.19 & 1.31 & 2 \\
& $\sigma_2$ & 0.01 & 0.03 & 1.55 & 3.63 & 0.00 & 1.22 & 1.20 & 1.27 & 2 \\
& $\sigma_3$ & 0.01 & 0.03 & 0.63 & 3.30 & 0.00 & 1.19 & 1.15 & 1.21 & 7 \\
& $\sigma_4$ & 0.01 & 0.09 & 2.15 & 4.63 & 0.82 & 1.29 & 1.23 & 1.37 & 0 \\
\midrule
40 & $\sigma_1$ & 0.03 & 0.58 & 1.30 & 2.34 & 0.20 & 1.20 & 1.16 & 1.25 & 0 \\
& $\sigma_2$ & 0.03 & 0.39 & 1.17 & 2.14 & 0.24 & 1.16 & 1.13 & 1.19 & 0 \\
& $\sigma_3$ & 0.02 & 0.18 & 0.89 & 1.87 & 0.00 & 1.18 & 1.13 & 1.27 & 1 \\
& $\sigma_4$ & 0.07 & 1.89 & 2.61 & 4.03 & 0.45 & 1.19 & 1.15 & 1.22 & 0 \\
\midrule
60 & $\sigma_1$ & 0.14 & 8.34 & 0.91 & 2.03 & 0.23 & 1.17 & 1.12 & 1.21 & 0 \\
& $\sigma_2$ & 0.12 & 1.62 & 0.98 & 1.90 & 0.34 & 1.13 & 1.10 & 1.15 & 0 \\
& $\sigma_3$ & 0.05 & 0.44 & 0.49 & 1.09 & 0.00 & 1.16 & 1.13 & 1.20 & 1 \\
& $\sigma_4$ & 0.39 & 30.41 & 2.57 & 3.97 & 0.94 & 1.14 & 1.12 & 1.16 & 0 \\
\midrule
80 & $\sigma_1$ & 0.41 & 7.96 & 0.74 & 1.88 & 0.28 & 1.14 & 1.11 & 1.17 & 0 \\
& $\sigma_2$ & 0.36 & 24.06 & 0.82 & 1.40 & 0.07 & 1.11 & 1.09 & 1.13 & 0 \\
& $\sigma_3$ & 0.19 & 1.89 & 0.47 & 0.85 & 0.17 & 1.15 & 1.12 & 1.19 & 0 \\
& $\sigma_4$ & 0.88 & \texttt{limit} & 5.78 & 10.77 & 2.20 & 1.11 & 1.10 & 1.12 & 0 \\
\midrule
100 & $\sigma_1$ & 0.73 & 8.76 & 0.46 & 0.82 & 0.06 & 1.13 & 1.11 & 1.14 & 0 \\
& $\sigma_2$ & 0.45 & 4.03 & 0.41 & 0.75 & 0.14 & 1.14 & 1.11 & 1.16 & 0 \\
& $\sigma_3$ & 0.32 & 0.81 & 0.17 & 0.68 & 0.00 & 1.15 & 1.12 & 1.20 & 2 \\
& $\sigma_4$ & 1.61 & \texttt{limit} & 4.44 & 7.66 & 1.34 & 1.09 & 1.08 & 1.10 & 0 \\
\bottomrule
\end{tabular}
\caption{Results for $\cgub$ and $\cglb$ with $C = 10$}
\label{tab:res10}
\end{table}

In Table~\ref{tab:res10},  it can be seen that, with $C=10$, the computation of $\cglb$
is {\em fast,\/} with average CPU times  less than 1 second in almost all the
cases (i.e., with any number of jobs).
The $\sigma_4$ instances are the most time demanding, with the only average
computation time above $1$~second. This is due to the fact that a larger set
of columns is usually generated on such instances.
The computation of $\cgub$ is, as expected, the heaviest part of the procedure, with 
larger CPU times. However,  only in cases $n=80,100$ and $\sigma=\sigma_4$ the CPLEX time limit
is reached. Again $\sigma_4$ instances were the most CPU time demanding, because
of the larger set of columns to be handled.
The certified solution quality was very good, with an average optimality gap
usually below $1.5\%$, and only one case ($n=80$, $\sigma=\sigma_4$) above
$5\%$.

From Table~\ref{tab:res10}, it can be easily seen that $\cglb$ performances are
much better than \texttt{PR} in every combination, ranging from an
average $9 \%$ gain when $n = 100$ and $\sigma = \sigma_4$ to an
average $29 \%$
when $n = 20$ and $\sigma = \sigma_4$. These values also suggest that
$\texttt{PR}$ performs better for large $n$;in fact, when a
high number of batches are required in the feasible solutions, the
usually weak parallel machine relaxation of \texttt{PR} becomes slightly
stronger.

\begin{table}[tbp]
\centering
\begin{tabular}{ll|rr|rrr|rrr|c}
\toprule
\multicolumn{2}{c}{Param} & \multicolumn{2}{c}{Times (s)} & \multicolumn{3}{c}{Gap (\%)} & \multicolumn{3}{c}{$\cglb / \texttt{PR}$} \\
\cmidrule(r){1-2}
\cmidrule(lr){3-4}
\cmidrule(lr){5-7}
\cmidrule(lr){8-10}
$n$ & $\sigma$ & \cglb & \cgub & avg & worst & best & avg & min & max & opt \\
\midrule
20 & $\sigma_1$ & 0.02 & 0.07 & 1.03 & 3.69 & 0.00 & 1.46 & 1.36 & 1.66 & 3 \\
& $\sigma_2$ & 0.02 & 0.08 & 1.28 & 3.77 & 0.00 & 1.39 & 1.29 & 1.49 & 5 \\
& $\sigma_3$ & 0.01 & 0.05 & 1.10 & 5.46 & 0.00 & 1.35 & 1.30 & 1.40 & 4 \\
& $\sigma_4$ & 0.02 & 0.09 & 0.00 & 0.00 & 0.00 & 1.81 & 1.62 & 2.11 & 10 \\
\midrule
40 & $\sigma_1$ & 0.21 & 3.11 & 3.13 & 6.43 & 0.21 & 1.30 & 1.23 & 1.39 & 0 \\
& $\sigma_2$ & 0.18 & 1.72 & 3.83 & 5.62 & 2.33 & 1.27 & 1.21 & 1.33 & 0 \\
& $\sigma_3$ & 0.10 & 2.95 & 4.37 & 6.29 & 2.98 & 1.20 & 1.17 & 1.26 & 0 \\
& $\sigma_4$ & 0.71 & 2.19 & 1.18 & 5.13 & 0.07 & 1.51 & 1.41 & 1.60 & 0 \\
\midrule
60 & $\sigma_1$ & 0.77 & \texttt{limit} & 6.78 & 10.27 & 3.51 & 1.22 & 1.18 & 1.27 & 0 \\
& $\sigma_2$ & 0.72 & 37.44 & 5.21 & 8.46 & 1.88 & 1.20 & 1.17 & 1.22 & 0 \\
& $\sigma_3$ & 0.42 & 30.35 & 3.74 & 5.79 & 1.44 & 1.15 & 1.13 & 1.18 & 0 \\
& $\sigma_4$ & 2.38 & 10.59 & 2.28 & 4.24 & 0.05 & 1.36 & 1.31 & 1.41 & 0 \\
\midrule
80 & $\sigma_1$ & 1.68 & \texttt{limit} & 11.05 & 17.40 & 3.40 & 1.21 & 1.16 & 1.23 & 0 \\
& $\sigma_2$ & 1.41 & \texttt{limit} & 11.68 & 41.32 & 2.65 & 1.16 & 1.13 & 1.20 & 0 \\
& $\sigma_3$ & 0.88 & \texttt{limit} & 7.94 & 12.24 & 3.53 & 1.12 & 1.11 & 1.13 & 0 \\
& $\sigma_4$ & 4.75 & \texttt{limit} & 7.01 & 18.67 & 2.41 & 1.29 & 1.27 & 1.32 & 0 \\
\midrule
100 & $\sigma_1$ & 3.27 & \texttt{limit} & 15.43 & 18.54 & 12.12 & 1.16 & 1.13 & 1.19 & 0 \\
& $\sigma_2$ & 2.87 & \texttt{limit} & 9.23 & 11.46 & 3.08 & 1.13 & 1.12 & 1.14 & 0 \\
& $\sigma_3$ & 1.79 & \texttt{limit} & 11.66 & 15.45 & 8.87 & 1.10 & 1.09 & 1.11 & 0 \\
& $\sigma_4$ & 8.85 & \texttt{limit} & 6.38 & 9.54 & 3.32 & 1.26 & 1.23 & 1.29 & 0 \\
\bottomrule
\end{tabular}
\caption{Results for $\cgub$ and $\cglb$ with $C = 30$}
\label{tab:res30}
\end{table}

From Table~\ref{tab:res30}, it can be noticed that CPU times for $\cglb$
increase; this is expected, since a larger number of possible batches
are generated with an increased capacity. The larger reduced master
problems obviously affect also the computation of $\cgub$, which reaches
the time limit in all the cases for $n=80,100$.
The average optimality gaps worsen, but the worse increase is not found
on $\sigma_4$ instances; instead, it affects more heavily $\sigma_1$ instances,
especially for large $n$.


Overall, increasing capacity also increments the distance between the
two lower bounds $\cglb$ and \texttt{PR}; $\cglb$ performs better in
every combination, ranging from an average $10 \%$ gain when $C = 30$,
$n = 100$, and $\sigma = \sigma_3$ to an average $81 \%$ when $C =
30$, $n = 20$, and $\sigma = \sigma_4$.  This is reasonable, since
\texttt{PR} is based on a preemptive relaxation to $C$ parallel
machines and allowing to split jobs on more machines weakens the
relaxation.

\begin{table}[tbp]
\centering
\begin{tabular}{ll|rr|rrr|rrr|c}
\toprule
\multicolumn{2}{c}{Param} & \multicolumn{2}{c}{Times (s)} & \multicolumn{3}{c}{Gap (\%)} & \multicolumn{3}{c}{$\cglb / \texttt{PR}$} \\
\cmidrule(r){1-2}
\cmidrule(lr){3-4}
\cmidrule(lr){5-7}
\cmidrule(lr){8-10}
$n$ & $\sigma$ & \cglb & \cgub & avg & worst & best & avg & min & max & opt \\
\midrule
20 & $\sigma_1$ & 0.02 & 0.09 & 0.00 & 0.00 & 0.00 & 1.70 & 1.48 & 1.94 & 10 \\
& $\sigma_2$ & 0.02 & 0.09 & 0.10 & 0.38 & 0.00 & 1.65 & 1.53 & 1.77 & 7 \\
& $\sigma_3$ & 0.02 & 0.08 & 0.28 & 2.66 & 0.00 & 1.52 & 1.35 & 1.62 & 7 \\
& $\sigma_4$ & 0.02 & 0.10 & 0.00 & 0.00 & 0.00 & 2.13 & 1.90 & 2.39 & 10 \\
\midrule
40 & $\sigma_1$ & 0.35 & 1.29 & 1.58 & 3.58 & 0.00 & 1.44 & 1.38 & 1.51 & 2 \\
& $\sigma_2$ & 0.44 & 1.31 & 1.78 & 3.16 & 0.35 & 1.41 & 1.34 & 1.54 & 0 \\
& $\sigma_3$ & 0.30 & 1.09 & 2.24 & 4.82 & 0.00 & 1.34 & 1.29 & 1.40 & 1 \\
& $\sigma_4$ & 0.75 & 2.45 & 0.11 & 0.81 & 0.00 & 1.80 & 1.67 & 1.87 & 8 \\
\midrule
60 & $\sigma_1$ & 1.39 & 23.59 & 4.47 & 6.64 & 1.24 & 1.38 & 1.30 & 1.44 & 0 \\
& $\sigma_2$ & 1.45 & 8.38 & 2.54 & 5.54 & 0.08 & 1.32 & 1.29 & 1.37 & 0 \\
& $\sigma_3$ & 0.76 & 13.37 & 4.30 & 6.65 & 1.81 & 1.23 & 1.20 & 1.27 & 0 \\
& $\sigma_4$ & 6.72 & 12.56 & 1.25 & 2.73 & 0.00 & 1.52 & 1.43 & 1.64 & 3 \\
\midrule
80 & $\sigma_1$ & 3.43 & \texttt{limit} & 7.40 & 10.26 & 3.17 & 1.32 & 1.27 & 1.38 & 0 \\
& $\sigma_2$ & 3.71 & \texttt{limit} & 3.95 & 4.64 & 2.63 & 1.26 & 1.24 & 1.27 & 0 \\
& $\sigma_3$ & 1.93 & \texttt{limit} & 4.94 & 10.32 & 2.00 & 1.19 & 1.17 & 1.20 & 0 \\
& $\sigma_4$ & 18.34 & \texttt{limit} & 1.90 & 4.85 & 0.11 & 1.46 & 1.42 & 1.51 & 0 \\
\midrule
100 & $\sigma_1$ & 7.09 & \texttt{limit} & 15.66 & 58.23 & 8.14 & 1.23 & 1.19 & 1.27 & 0 \\
& $\sigma_2$ & 7.04 & \texttt{limit} & 7.19 & 10.33 & 4.27 & 1.21 & 1.19 & 1.23 & 0 \\
& $\sigma_3$ & 3.50 & \texttt{limit} & 7.86 & 11.47 & 3.51 & 1.15 & 1.14 & 1.17 & 0 \\
& $\sigma_4$ & 38.65 & \texttt{limit} & 3.26 & 5.70 & 0.84 & 1.39 & 1.35 & 1.47 & 0 \\
\bottomrule
\end{tabular}
\caption{Results for $\cgub$ and $\cglb$ with $C = 50$}
\label{tab:res50}
\end{table}

Table~\ref{tab:res50} shows the results of the tests with capacity $C=50$ that
confirm the impact of $C$. The instances belonging to class $\sigma_4$ are still
the most computationally demanding, both for lower bounding and heuristic
solution. Instances with $\sigma=\sigma_1$ are the worse in terms
of solution quality --- with the exception of small $20$ job instances ---
but, curiously, the gap lowers on $n=80$ instances
when passing from $C=30$ to $C=50$. The worst average gap is reached with
$n=100$ and $\sigma=\sigma_1$ ($15.66\%)$.
Also, \texttt{PR} worsens considerably with respect to $\cglb$.

\begin{table}[tbp]
\centering
\begin{tabular}{ll|rr|rr}
\toprule
\multicolumn{2}{c}{Param} & \multicolumn{2}{c}{Times (s)} & \multicolumn{2}{c}{$\text{Heur}/ \texttt{PR}$} \\
\cmidrule(r){1-2}
\cmidrule(lr){3-4}
\cmidrule(lr){5-6}
$n$ & $\sigma$ & \texttt{HMMAS} & \cgub & \texttt{HMMAS} & \cgub \\
\midrule
20 & $\sigma_1$ & 1.29 & 0.04 & 1.25 & 1.27 \\
& $\sigma_2$ & 1.41 & 0.03 & 1.25 & 1.24 \\
& $\sigma_3$ & 1.56 & 0.03 & 1.21 & 1.20 \\
& $\sigma_4$ & 0.98 & 0.09 & 1.28 & 1.31 \\
\midrule
40 & $\sigma_1$ & 5.55 & 0.58 & 1.19 & 1.21 \\
& $\sigma_2$ & 5.93 & 0.39 & 1.19 & 1.18 \\
& $\sigma_3$ & 6.32 & 0.18 & 1.18 & 1.19 \\
& $\sigma_4$ & 3.68 & 1.89 & 1.20 & 1.22 \\
\midrule
60 & $\sigma_1$ & 20.61 & 8.34 & 1.17 & 1.18 \\
& $\sigma_2$ & 19.26 & 1.62 & 1.16 & 1.14 \\
& $\sigma_3$ & 23.41 & 0.44 & 1.18 & 1.17 \\
& $\sigma_4$ & 10.90 & 30.41 & 1.18 & 1.17 \\
\midrule
80 & $\sigma_1$ & 57.54 & 7.96 & 1.16 & 1.15 \\
& $\sigma_2$ & 55.72 & 24.06 & 1.16 & 1.12 \\
& $\sigma_3$ & 62.91 & 1.89 & 1.16 & 1.15 \\
& $\sigma_4$ & 28.87 & 60.91 & 1.16 & 1.18 \\
\midrule
100 & $\sigma_1$ & 109.72 & 8.76 & 1.16 & 1.13 \\
& $\sigma_2$ & 106.39 & 4.03 & 1.15 & 1.14 \\
& $\sigma_3$ & 135.90 & 0.81 & 1.15 & 1.16 \\
& $\sigma_4$ & 53.58 & 61.66 & 1.15 & 1.14 \\
\bottomrule
\end{tabular}
\caption{Comparison between \texttt{HMMAS} and $\cgub$ algorithms}
\label{tab:ant}
\end{table}

An attempt to compare our upper bound $\cgub$ to the
hybrid ant system (\texttt{HMMAS})
developed by~\cite{IranianiFormiche} has  been done by generating the random instances in the
same way they did, and evaluating the execution times with extreme
care considering also the different operational environment, as  neither their algorithm nor to their instances were available.
Also, for the comparison, only capacity $C =10$ was used, as in~\cite{IranianiFormiche} only such value has been used.

As can be seen in Table~\ref{tab:ant}, the results show that the
performance of $\cgub$, evaluated against Uzsoy's lower bound
\texttt{PR}, seems to be very similar to that of \texttt{HMMAS}.
It is not possible to explicitly compare to results of~\cite{IranianiFormiche}
but, since the results appear to be very close, it can be speculated that
the two algorithms could give similar results for the upper bound,
when they are run on the same instance set.
With the same care, it can be observed that the CPU times of $\cgub$ seem to be
much smaller than those of \texttt{HMMAS}, even taking into account the
different processors. The notable exception is for large $\sigma_4$
instances that generate large reduced master problems.


It must be stressed, however, that, as the instances
of~\cite{IranianiFormiche} were not available, the optimality gap of
their results against a strong lower bound is unknown. Thus, even if
the results seem to suggest that the upper bounds are comparable, the
algorithm quality cannot be directly benchmarked.

Eventually, the quality of $\cgub$ has been compared to the quality of the
heuristic solution reached by CPLEX (\texttt{CPLEX-UB}) after $300$
seconds of computation
using model~\eqref{mip-obj}--\eqref{mip-vars}. CPLEX optimality gap
is most of the times well above $90\%$ because the lower bound is zero
or almost zero. Anyway, using the proposed stronger lower bound,  
 a more realistic optimality gap can be computed for CPLEX as
\[\frac{\texttt{CPLEX-UB}-\cglb}{\tt UB^\ast}\cdot100\%,
\qquad\text{with ${\tt UB^\ast}=\min\{\texttt{CPLEX-UB},\cgub\}$}.\]
The gap for $\cgub$ is recomputed as $(\cgub-\cglb)/{\tt UB^\ast}\cdot100\%$ for
uniformity.

The comparison is reported in Table~\ref{tab:cplex}, again in terms
of average, worst and best gap. The column $\#\text{win}$ counts the number
of instances out of ten for which each algorithm achieves the best solution.
In case of a draw,  a ``win'' is counted for both, so the two columns
can sum to more than $10$.
Instances with $n=20$, $40$, $60$, $80$ and $C=10$ have been tested. CPLEX ran for
the full $300$ seconds on all the instances, without proving optimality
for any of them.
$\cgub$ ran with the same $60$ seconds time limit as in
Table~\ref{tab:res10}. Basically, except for the small $n=20$ instances,
CPLEX solution is consistently worse than $\cgub$.

\begin{table}[tbp]
\centering
\begin{tabular}{ll|rrrc|rrrc}
\toprule
\multicolumn{2}{c}{Param} & \multicolumn{3}{c}{\texttt{CPLEX-UB} Gap (\%)}
& \multicolumn{1}{c}{} & \multicolumn{3}{c}{$\cgub$ Gap (\%)} \\
\cmidrule(r){1-2}
\cmidrule(lr){3-5}
\cmidrule(lr){7-9}
$n$ & $\sigma$ & avg & worst & best & \#win &
avg & worst & best & \#win \\
\midrule
20 & $\sigma_1$ & 1.15 & 3.73 & 0.00 & 8 & 1.31 & 3.27 & 0.00 & 6 \\
& $\sigma_2$ & 1.17 & 3.58 & 0.00 & 7 & 1.56 & 3.73 & 0.00 & 7 \\
& $\sigma_3$ & 0.69 & 3.30 & 0.00 & 8 & 0.63 & 3.30 & 0.00 & 10 \\
& $\sigma_4$ & 3.51 & 7.14 & 0.82 & 4 & 2.15 & 4.63 & 0.82 & 9 \\
\midrule
40 & $\sigma_1$ & 9.60 & 14.63 & 5.97 & 0 & 1.30 & 2.34 & 0.20 & 10 \\
& $\sigma_2$ & 9.34 & 16.02 & 5.43 & 0 & 1.17 & 2.14 & 0.24 & 10 \\
& $\sigma_3$ & 4.41 & 8.13 & 2.34 & 0 & 0.89 & 1.87 & 0.00 & 10 \\
& $\sigma_4$ & 12.03 & 18.16 & 9.09 & 0 & 2.61 & 4.03 & 0.45 & 10 \\
\midrule
60 & $\sigma_1$ & 49.28 & 77.33 & 38.61 & 0 & 0.91 & 2.03 & 0.23 & 10 \\
& $\sigma_2$ & 48.86 & 59.69 & 33.98 & 0 & 0.98 & 1.90 & 0.34 & 10 \\
& $\sigma_3$ & 33.78 & 41.74 & 22.84 & 0 & 0.49 & 1.09 & 0.00 & 10 \\
& $\sigma_4$ & 45.86 & 66.77 & 24.59 & 0 & 2.57 & 3.97 & 0.94 & 10 \\
\midrule
80 & $\sigma_1$ & 73.77 & 88.55 & 59.86 & 0 & 0.74 & 1.88 & 0.28 & 10 \\
& $\sigma_2$ & 66.12 & 77.45 & 53.45 & 0 & 0.82 & 1.40 & 0.07 & 10 \\
& $\sigma_3$ & 52.79 & 73.82 & 38.40 & 0 & 0.47 & 0.85 & 0.17 & 10 \\
& $\sigma_4$ & 84.09 & 102.96 & 62.95 & 0 & 5.78 & 10.77 & 2.20 & 10 \\
\bottomrule
\end{tabular}
\caption{Comparison between \texttt{CPLEX-UB} (300 secs) and $\cgub$.}
\label{tab:cplex}
\end{table}

\subsection{Parallel machines}

With the same data, the $Pm|\pbatch|\sum C_j$ problem has been solved with
$m=2,3$ and~$5$ machines. The testing has been limited to the case $C=10$.
The time limit for the branch and bound phase of the heuristic was
raised to $180$ seconds.  The results are reported in tables~%
\ref{tab:res10m2}, \ref{tab:res10m3} and~\ref{tab:res10m5}.
Apparently, increasing the number of machines has a very mild impact
on the CPU time for computing the lower bound.  The growth of the
computational cost is much higher for the branch and bound phase, but
with a certain variability on the four classes of instances, with
classes $\sigma_1$ and $\sigma_4$ exhibiting the largest
growth. Again, class $\sigma_4$ broke the time limit in all the instances.  The
quality of the solution, as measured by the percentage gap, does not
suffer seriously, except for case $n=100$, $m=5$, class $\sigma_4$.  The
worst average of $14.09\%$  is caused by one single instance with
a very large gap of $81.62\%$; if a larger but still acceptable
time limit of $300$ seconds is allowed, the average gap for this class
lowers to $4.63\%$ (max gap $12.56\%$).

Uzsoy's bound $\texttt{PR}$ is easily extended to the parallel machines case
allowing a relaxation to $mC$ parallel machines.
Tables~\ref{tab:res10m2}, \ref{tab:res10m3} and~\ref{tab:res10m5}
also compare $\cglb$ with 
{\tt PR} extended to the parallel machines case. The
ratio between the two bounds is apparently unaffected by the growth of
$m$.

\section{Final remarks}
\label{concl}
In this paper, column generation techniques for solving $1|\pbatch|\sum C_j$
problem has been explored, generalizing such techniques to problems with parallel machines.
The exponential size model~\eqref{cg:obj}--\eqref{cg:vars}, handled by
means of column generation, allows to find --- to the authors' knowledge ---
the tightest known lower bound for $1|\pbatch|\sum C_j$. Embedded in a simple
price-and-branch approach, it achieves high-quality solutions for instances
up to $100$ jobs in size, with certified optimality gaps.
The model relies on Property~\ref{prop1} in order to express the linear
objective function by means of ``positional'' coefficients.
 Property~\ref{prop2} is crucial in order to develop
an efficient pricing procedure.
The approach is flexible enough to be extended to problems with
parallel machines with a very limited effort 
while it is probably not so simple to extend the it to weighted $\sum w_jC_j$
objectives.

Having a column-generation based lower bound would naturally lead to
searching for a branch-and-price exact approach; some problems have
still to be tackled in this direction.
Applying the classical branching scheme of~\citep{RyanFoster}
makes the pricing procedure useless at most branching nodes,
since the knapsack-like
problems must take disjunctions into account. 

Preliminary experiments have been run
with a branching scheme that allows to keep the pricing problem structure.
Similarly to what is known to happen with other $\sum C_j$ problems,
a narrow optimality gap is obtained at the root node,
but closing such gap is quite difficult.
In the preliminary experiments, this
seems to happen because of the large number of equivalent optimal solutions
that can potentially be generated playing on the different equivalent packing
of jobs into batches. This leads to a large number of optimal branches that
cannot be fathomed by bounding --- at least not in early stages of the search.
Some strong dominance criterion, suitable to break or prune
such equivalencies, might then be needed before a
branch-and-price is usable for this type of parallel batching problems.
This is subject of ongoing research.

\begin{table}[htbp]
\centering
\begin{tabular}{ll|rr|rrr|rrr|c}
\toprule
\multicolumn{2}{c}{Param} & \multicolumn{2}{c}{Times (s)} & \multicolumn{3}{c}{Gap (\%)} & \multicolumn{3}{c}{$\cglb / \texttt{PR}$} \\
\cmidrule(r){1-2}
\cmidrule(lr){3-4}
\cmidrule(lr){5-7}
\cmidrule(lr){8-10}
$n$ & $\sigma$ & \cglb & \cgub & avg & worst & best & avg & min & max & opt \\
\midrule
20 & $\sigma_1$ & 0.03 & 0.12 & 0.78 & 2.00 & 0.00 & 1.24 & 1.18 & 1.29 & 2 \\
& $\sigma_2$ & 0.03 & 0.09 & 1.43 & 4.78 & 0.00 & 1.21 & 1.19 & 1.25 & 2 \\
& $\sigma_3$ & 0.02 & 0.07 & 0.57 & 3.04 & 0.00 & 1.18 & 1.14 & 1.19 & 7 \\
& $\sigma_4$ & 0.04 & 0.22 & 2.08 & 4.11 & 0.00 & 1.28 & 1.22 & 1.35 & 1 \\
\midrule
40 & $\sigma_1$ & 0.06 & 0.82 & 1.14 & 2.28 & 0.01 & 1.19 & 1.16 & 1.24 & 0 \\
& $\sigma_2$ & 0.05 & 0.90 & 1.08 & 1.73 & 0.23 & 1.16 & 1.13 & 1.18 & 0 \\
& $\sigma_3$ & 0.03 & 0.22 & 0.75 & 1.39 & 0.00 & 1.18 & 1.13 & 1.26 & 1 \\
& $\sigma_4$ & 0.13 & 5.01 & 2.07 & 3.73 & 0.37 & 1.18 & 1.15 & 1.22 & 0 \\
\midrule
60 & $\sigma_1$ & 0.25 & 10.22 & 0.83 & 1.64 & 0.23 & 1.16 & 1.12 & 1.20 & 0 \\
& $\sigma_2$ & 0.21 & 1.72 & 0.90 & 1.56 & 0.32 & 1.12 & 1.10 & 1.15 & 0 \\
& $\sigma_3$ & 0.13 & 0.60 & 0.38 & 1.00 & 0.01 & 1.16 & 1.13 & 1.20 & 0 \\
& $\sigma_4$ & 0.49 & 76.71 & 2.35 & 3.72 & 1.69 & 1.14 & 1.12 & 1.16 & 0 \\
\midrule
80 & $\sigma_1$ & 0.62 & 45.96 & 0.76 & 1.94 & 0.21 & 1.14 & 1.10 & 1.17 & 0 \\
& $\sigma_2$ & 0.54 & 40.97 & 0.71 & 1.00 & 0.08 & 1.11 & 1.09 & 1.13 & 0 \\
& $\sigma_3$ & 0.35 & 3.29 & 0.46 & 0.76 & 0.18 & 1.14 & 1.12 & 1.18 & 0 \\
& $\sigma_4$ & 1.30 & \texttt{limit} & 5.65 & 10.31 & 1.33 & 1.11 & 1.10 & 1.12 & 0 \\
\midrule
100 & $\sigma_1$ & 1.08 & 19.16 & 0.44 & 0.77 & 0.07 & 1.13 & 1.11 & 1.14 & 0 \\
& $\sigma_2$ & 0.79 & 3.66 & 0.41 & 0.73 & 0.12 & 1.13 & 1.11 & 1.16 & 0 \\
& $\sigma_3$ & 0.59 & 1.63 & 0.21 & 0.66 & 0.00 & 1.15 & 1.11 & 1.20 & 2 \\
& $\sigma_4$ & 2.25 & \texttt{limit} & 7.19 & 26.27 & 1.08 & 1.09 & 1.08 & 1.10 & 0 \\
\bottomrule
\end{tabular}
\caption{Results for $\cgub$ and $\cglb$ with $C = 10$ and 2 parallel machines}
\label{tab:res10m2}
\end{table}

\begin{table}[htbp]
\centering
\begin{tabular}{ll|rr|rrr|rrr|c}
\toprule
\multicolumn{2}{c}{Param} & \multicolumn{2}{c}{Times (s)} & \multicolumn{3}{c}{Gap (\%)} & \multicolumn{3}{c}{$\cglb / \texttt{PR}$} \\
\cmidrule(r){1-2}
\cmidrule(lr){3-4}
\cmidrule(lr){5-7}
\cmidrule(lr){8-10}
$n$ & $\sigma$ & \cglb & \cgub & avg & worst & best & avg & min & max & opt \\
\midrule
20 & $\sigma_1$ & 0.03 & 0.13 & 0.52 & 1.29 & 0.00 & 1.24 & 1.18 & 1.29 & 3 \\
& $\sigma_2$ & 0.03 & 0.10 & 1.15 & 2.90 & 0.00 & 1.21 & 1.19 & 1.25 & 2 \\
& $\sigma_3$ & 0.03 & 0.08 & 0.45 & 2.49 & 0.00 & 1.18 & 1.14 & 1.19 & 7 \\
& $\sigma_4$ & 0.05 & 0.30 & 1.77 & 3.61 & 0.00 & 1.30 & 1.23 & 1.36 & 1 \\
\midrule
40 & $\sigma_1$ & 0.12 & 1.68 & 1.11 & 2.02 & 0.00 & 1.19 & 1.16 & 1.24 & 1 \\
& $\sigma_2$ & 0.09 & 1.00 & 0.88 & 1.38 & 0.00 & 1.16 & 1.13 & 1.18 & 1 \\
& $\sigma_3$ & 0.06 & 0.43 & 0.80 & 1.64 & 0.00 & 1.17 & 1.12 & 1.25 & 1 \\
& $\sigma_4$ & 0.19 & 4.81 & 1.76 & 3.37 & 0.19 & 1.19 & 1.16 & 1.22 & 0 \\
\midrule
60 & $\sigma_1$ & 0.35 & 3.89 & 0.73 & 1.54 & 0.24 & 1.16 & 1.12 & 1.20 & 0 \\
& $\sigma_2$ & 0.32 & 3.47 & 0.86 & 1.50 & 0.29 & 1.12 & 1.10 & 1.15 & 0 \\
& $\sigma_3$ & 0.24 & 1.09 & 0.38 & 1.12 & 0.00 & 1.15 & 1.13 & 1.19 & 1 \\
& $\sigma_4$ & 0.62 & 76.82 & 2.00 & 3.92 & 1.01 & 1.15 & 1.12 & 1.16 & 0 \\
\midrule
80 & $\sigma_1$ & 1.04 & 16.93 & 0.62 & 1.57 & 0.18 & 1.14 & 1.10 & 1.16 & 0 \\
& $\sigma_2$ & 0.75 & 47.57 & 0.71 & 1.12 & 0.08 & 1.11 & 1.09 & 1.12 & 0 \\
& $\sigma_3$ & 0.58 & 3.99 & 0.47 & 0.72 & 0.19 & 1.14 & 1.12 & 1.18 & 0 \\
& $\sigma_4$ & 1.49 & \texttt{limit} & 5.23 & 13.93 & 1.20 & 1.12 & 1.10 & 1.13 & 0 \\
\midrule
100 & $\sigma_1$ & 1.99 & 20.33 & 0.44 & 0.73 & 0.07 & 1.12 & 1.11 & 1.14 & 0 \\
& $\sigma_2$ & 1.12 & 4.92 & 0.39 & 0.69 & 0.12 & 1.13 & 1.10 & 1.16 & 0 \\
& $\sigma_3$ & 1.18 & 2.83 & 0.17 & 0.65 & 0.00 & 1.15 & 1.11 & 1.20 & 2 \\
& $\sigma_4$ & 3.09 & \texttt{limit} & 5.44 & 19.24 & 1.58 & 1.09 & 1.08 & 1.10 & 0 \\
\bottomrule
\end{tabular}
\caption{Results for $\cgub$ and $\cglb$ with $C = 10$ and 3 parallel machines}
\label{tab:res10m3}
\end{table}

\begin{table}[htbp]
\centering
\begin{tabular}{ll|rr|rrr|rrr|c}
\toprule
\multicolumn{2}{c}{Param} & \multicolumn{2}{c}{Times (s)} & \multicolumn{3}{c}{Gap (\%)} & \multicolumn{3}{c}{$\cglb / \texttt{PR}$} \\
\cmidrule(r){1-2}
\cmidrule(lr){3-4}
\cmidrule(lr){5-7}
\cmidrule(lr){8-10}
$n$ & $\sigma$ & \cglb & \cgub & avg & worst & best & avg & min & max & opt \\
\midrule
20 & $\sigma_1$ & 0.05 & 0.20 & 0.47 & 1.95 & 0.00 & 1.27 & 1.18 & 1.30 & 2 \\
& $\sigma_2$ & 0.05 & 0.16 & 0.84 & 1.77 & 0.00 & 1.22 & 1.20 & 1.24 & 3 \\
& $\sigma_3$ & 0.05 & 0.10 & 0.30 & 1.57 & 0.00 & 1.19 & 1.15 & 1.20 & 7 \\
& $\sigma_4$ & 0.05 & 0.26 & 0.83 & 2.20 & 0.00 & 1.36 & 1.27 & 1.44 & 2 \\
\midrule
40 & $\sigma_1$ & 0.23 & 1.94 & 0.91 & 1.59 & 0.24 & 1.19 & 1.17 & 1.23 & 0 \\
& $\sigma_2$ & 0.18 & 1.62 & 0.71 & 1.25 & 0.00 & 1.16 & 1.14 & 1.19 & 1 \\
& $\sigma_3$ & 0.15 & 0.68 & 0.72 & 1.79 & 0.02 & 1.17 & 1.13 & 1.24 & 0 \\
& $\sigma_4$ & 0.27 & 6.02 & 1.28 & 2.43 & 0.32 & 1.21 & 1.17 & 1.26 & 0 \\
\midrule
60 & $\sigma_1$ & 0.53 & 5.34 & 0.70 & 1.25 & 0.25 & 1.16 & 1.13 & 1.20 & 0 \\
& $\sigma_2$ & 0.55 & 7.23 & 0.81 & 1.78 & 0.34 & 1.12 & 1.10 & 1.15 & 0 \\
& $\sigma_3$ & 0.45 & 2.17 & 0.36 & 0.89 & 0.00 & 1.15 & 1.13 & 1.19 & 1 \\
& $\sigma_4$ & 0.69 & 89.00 & 1.65 & 3.14 & 0.37 & 1.16 & 1.13 & 1.18 & 0 \\
\midrule
80 & $\sigma_1$ & 1.33 & 17.41 & 0.61 & 1.56 & 0.11 & 1.14 & 1.11 & 1.16 & 0 \\
& $\sigma_2$ & 0.93 & 38.02 & 0.62 & 0.99 & 0.01 & 1.11 & 1.09 & 1.12 & 0 \\
& $\sigma_3$ & 0.78 & 4.37 & 0.39 & 0.60 & 0.16 & 1.14 & 1.12 & 1.17 & 0 \\
& $\sigma_4$ & 1.53 & \texttt{limit} & 2.50 & 4.57 & 1.27 & 1.13 & 1.11 & 1.14 & 0 \\
\midrule
100 & $\sigma_1$ & 2.76 & 20.14 & 0.43 & 0.79 & 0.05 & 1.12 & 1.11 & 1.14 & 0 \\
& $\sigma_2$ & 1.41 & 6.34 & 0.38 & 0.61 & 0.12 & 1.13 & 1.10 & 1.15 & 0 \\
& $\sigma_3$ & 1.54 & 3.50 & 0.16 & 0.63 & 0.00 & 1.14 & 1.11 & 1.19 & 2 \\
& $\sigma_4$ & 3.22 & \texttt{limit} & 14.09 & 81.62 & 0.87 & 1.10 & 1.09 & 1.11 & 0 \\
\bottomrule
\end{tabular}
\caption{Results for $\cgub$ and $\cglb$ with $C = 10$ and 5 parallel machines}
\label{tab:res10m5}
\end{table}
\clearpage 

\bibliographystyle{elsarticle-num-names} 
\bibliography{biblio}

\end{document}